\documentclass[psamsfonts]{amsart}

\usepackage{amssymb,amsfonts,amscd}

\usepackage[colorlinks,linktocpage]{hyperref}
\hypersetup{linkcolor=[rgb]{0,0,0.715}}
\hypersetup{citecolor=[rgb]{0,0.715,0}}

\newtheorem{thm}{Theorem}[section]

\newtheorem*{thm1}{Theorem A}

\theoremstyle{definition}

\theoremstyle{remark}
\newtheorem{rem}[thm]{Remark}

\makeatletter
\let\c@equation\c@thm
\makeatother
\numberwithin{equation}{section}

\title[]{The Grushin hemisphere as a Ricci limit space\\ with curvature $\ge 1$}
\author[]{Jiayin Pan}
\address[Jiayin Pan]{Department of Mathematics, University of California, Santa Curz, CA, US.}
\email{jpan53@ucsc.edu}

\begin{document}
	
\begin{abstract}
	The Grushin sphere is an almost-Riemannian manifold that degenerates along its equator. We construct a sequence of Riemannian metrics on a sphere $S^{m+n}$ with $Ric\ge 1$ such that its Gromov-Hausdorff limit is the $n$-dimensional Grushin hemisphere.
\end{abstract}    	
	
\maketitle

A Ricci limit space is the Gromov-Hausdorff limit of a sequence of complete Riemannian $n$-manifolds with Ricci curvature uniformly bounded below. The structure of Ricci limit spaces is crucial in understanding Ricci curvature and has been studied extensively since the seminal works by Cheeger and Colding. In a previous joint work with Wei \cite{PanWei}, we have constructed a Ricci limit space, as the asymptotic cone of a complete Riemannian metric on $\mathbb{R}^{n+1}$ with $Ric>0$, such that its Hausdorff dimension exceeds its rectifiable dimension. The limit space is a halfplane $[0,\infty)\times \mathbb{R}$; its singular set is the boundary $\{0\}\times \mathbb{R}$ with Hausdorff dimension $1+\alpha$, where $\alpha>0$ can be any prior chosen number. It answered a longstanding open problem by Cheeger and Colding \cite{ChCo97}. In a joint work with Dai, Honda, and Wei \cite{DHPW}, among other results, we have given a more detailed description of these spaces: it is the metric completion of an incomplete weighted Riemannian metric defined on the open halfplane $(0,\infty)\times \mathbb{R}$. The metric and measure are
$$g=dx^2+x^{-2\alpha} dy^2,\quad \mathfrak{m}=c x^{\frac{n-1}{2}-\alpha}dxdy.$$

A subRiemannian manifold is a manifold endowed with a distribution and a fiberwise inner product on the distribution. The distribution specifies in which directions one can travel. A classical example is the Heisenberg $3$-group, which is a nilpotent group with topological dimension $3$ and Hausdorff dimension $4$. Another simple example is the Grushin plane (see \cite[Section 3.1]{Bel}). Its distribution is generated by the vector fields
$X=\partial_x$ and $Y=|x|^\alpha \partial_y$ on $\mathbb{R}^2$,
where $\alpha>0$. Setting $\{X,Y\}$ orthonormal defines a subRiemannian metric on the plane. Note that $Y$ only degenerates along the $y$-axis, so the distribution has maximal rank almost everywhere. A subRiemannian manifold with this property is called almost-Riemannian. Outside the $y$-axis, the Grushin plane becomes Riemannian with metric
$$g=dx^2+|x|^{-2\alpha} dy^2.$$
Its halfplane $[0,\infty)\times \mathbb{R}$ is a convex subset, and the boundary $\{0\}\times \mathbb{R}$ has Hausdorff dimension $1+\alpha$.

For readers, it is now clear that Pan-Wei's Ricci limit examples, as metric spaces, are isometric to the Grushin halfplanes. This surprising connection was pointed out to the author by Richard Montgomery when they met each other for the first time. Special thanks to him. It was a moment Ricci curvature met subRiemannian geometry.

In this paper, now we turn to another almost-Riemannian manifold: the Grushin sphere, which is the spherical analog of the Grushin plane \cite{BCGGJ}. The distribution of the $2$-dimensional Grushin sphere is generated by the vector fields
$X=(1,0)$ and $Y=(0,\tan\phi)$
written in the spherical coordinate $(\phi,\theta)\in(-\pi/2,\pi/2)\times[0,2\pi]$. Setting $\{X,Y\}$ orthonormal defines the subRiemannian metric. The equator $\{\phi=0\}$ has Hausdorff dimension $2$. The Grushin hemisphere $\{\phi\ge 0\}$ is convex and can be viewed as the metric completion of a warped product on the open hemisphere
$$(0,\pi/2]\times_h S^{1},\quad d\phi^2+(\tan \phi)^{-2} d\theta^2$$
The $n$-dimensional Grushin sphere are constructed similarly. Its equator has Hausdorff dimension $2(n-1)$. It seems to the author that we don't have a good notion for the $\alpha$-variants; using $h(\phi)=(\tan \phi)^{-\alpha}$ as the warping function produces singularities at the poles since $h'(\pi/2)\not=-1$ when $\alpha\not=1$. 

Inspired by the above-mentioned connection between Ricci curvature and subRiemannian geometry, we construct a sequence of Riemannian metrics on a sphere with $Ric\ge 1$ converging to the $n$-dimensional Grushin hemisphere.

\begin{thm1}
	Given an integer $n\ge 2$, there is a sequence of Riemannian metrics $g_i$ on $S^{m+n}$ with $Ric(g_i)\ge 1$, where $m$ is sufficiently large, such that $(S^{m+n},g_i)$ Gromov-Hausdorff converges to the $n$-dimensional Grushin hemisphere. 
\end{thm1}

As a Ricci limit space, the Grushin hemisphere carries a limit renormalized measure from the sequence. We emphasize that this limit renormalized measure is different from the measure induced by the almost-Riemannian metric.

It had been long believed that subRiemannian geometry cannot interact with Ricci limit spaces (or the RCD$(K,N)$ condition, or the even weaker CD$(K,N)$ condition). In fact, by a result of Juillet \cite{Ju}, any complete subRiemannian manifold with a distribution of non-maximal rank everywhere and a measure of smooth positive density does not satisfy the CD$(K,N)$ condition for any $K$ and $N$. By a recent result of Magnabosco and Rossi \cite{MR}, any complete $2$-dimensional almost-Riemannian manifold with a measure of smooth positive density does not satisfy the CD$(K,N)$ condition for any $K$ and $N$.

Now with the Grushin halfplanes and the Grushin hemispheres constructed as Ricci limit spaces, we establish a surprising connection.

Comparing Theorem A with Pan-Wei's construction \cite{PanWei}, one of the main differences here is the positive Ricci curvature lower bound. In \cite{PanWei}, compact Ricci limit spaces with large Hausdorff dimension were constructed, but they must have negative Ricci curvature somewhere. On a technical note, we remark that the fundamental group plays an essential role in Pan-Wei's construction, while in Theorem A we directly construct the metrics on a sphere.

\textit{Acknowledgements:} The author would like to thank Richard Montgomery for introducing him on the Grushin plane. The author would like to thank Nic Brody, Xianzhe Dai, Shouhei Honda, and Guofang Wei for helpful discussions. The author would like to thank Aaron Naber abd Luca Ruzzi for helpful comments.

\section{Construction}\label{sec_1}
Let $n\ge 2$ and $m$ to be determined later. We construct a family of Riemannian metrics $\{g_\lambda\}_{\lambda\ge 1}$ as doubly warped products:
$$M=[0,\pi/2]\times_{f_\lambda} S^{m} \times_{h_\lambda} S^{n-1},\quad g_\lambda=dr^2+f_\lambda(r)^2ds_m^2+h_\lambda(r)^2ds_{n-1}^2,$$
where $ds_k^2$ denotes the standard metric on the sphere $S^k$.
We use the warping functions as
$$f_\lambda(r)=\dfrac{\sin r}{(1+\lambda^2 \sin^2 r)^{1/4}},\quad h_\lambda(r)=\left( \dfrac{1}{\lambda^2}+\tan^2 r \right)^{-1/2}.$$
$f_\lambda$ and $h_\lambda$ satisfy
$$f_\lambda(0)=0, \quad f'_\lambda(0)=1,\quad f^{(\text{even})}_\lambda(0)=0,\quad f_\lambda(\pi/2)>0,\quad f^{(\text{odd})}_\lambda(\pi/2)=0;$$
$$h_\lambda(0)>0, \quad h^{(\text{odd})}_\lambda(0)=0,\quad h_\lambda(\pi/2)=0,\quad h'_\lambda(\pi/2)=-1,\quad h^{(\text{even})}_\lambda(\pi/2)=0.$$
Therefore, $g_\lambda$ defines a smooth Riemannian metric on $M$. Topologically, $M$ is a quotient
$$M=[0,\pi/2]\times S^m \times S^{n-1}/\sim,$$
where $\sim$ is given by
$$(0,x,y)\sim(0,x',y),\quad (\pi/2,x,y)\sim (\pi/2,x,y')$$
for all $x,x'\in S^m$ and all $y,y'\in S^{n-1}$. $M$ is diffeomorphic to the sphere $S^{m+n}$.
 
Let $H=\partial_r$, $U$ a unit vector tangent to $S^{m}$, $V$ a unit vector tangent to $S^{n-1}$. Then $(M,g_\lambda)$ has Ricci curvature 
\begin{align*}
	Ric(H,H)=& -m\dfrac{f_\lambda''}{f_\lambda}-(n-1)\dfrac{h_\lambda''}{h_\lambda}, \\
	Ric(U,U)=& -\dfrac{f''_\lambda}{f_\lambda}+(m-1)\dfrac{1-f'^2_\lambda}{f^2_\lambda}-(n-1)\dfrac{f_\lambda'h_\lambda'}{f_\lambda h_\lambda}, \\
	Ric(V,V)=& -\dfrac{h''_\lambda}{h_\lambda}+(n-2)\dfrac{1-h'^2_\lambda}{h^2_\lambda}-m\dfrac{f'_\lambda h'_\lambda}{f_\lambda h_\lambda}.
\end{align*}
We provide detailed calculations here. For convenience, we write
$$A=\lambda^2\sin^2r+\cos^2r,\quad B=\lambda^2\sin^2r+1.$$

By direct calculation and simplification, we have
$$f'_\lambda=\dfrac{(\cos r)(B+1)}{2B^{5/4}} \in [0,1].$$
$$-\dfrac{f''_\lambda}{f_\lambda}=\dfrac{\lambda^4\sin^4 r+ \lambda^4\sin^2 r+6\lambda^2+4}{4B^2}\ge 1.$$
$$-\dfrac{f'_\lambda h'_\lambda}{f_\lambda h_\lambda}=\dfrac{\lambda^2(B+1)}{2AB}\ge \dfrac{1}{2}.$$
The above three imply that $Ric(U,U)\ge 1$. 

Next, we check $Ric(H,H)$. We have
$$-\dfrac{h''_\lambda}{h_\lambda}=\dfrac{\lambda^2}{A^2}\left( (-2\lambda^2+2)\sin^2r+1 \right)\ge -\dfrac{2\lambda^4\sin^2r}{A^2}+1.$$
\begin{align*}
	&Ric(H,H)=-m\dfrac{f_\lambda''}{f_\lambda}-(n-1)\dfrac{h_\lambda''}{h_\lambda}\\
	\ge&\ m\cdot \dfrac{\lambda^4\sin^4 r+ \lambda^4\sin^2 r+6\lambda^2+4}{4B^2}-(n-1)\dfrac{2\lambda^4\sin^2r}{A^2}+1\\
	=&\ \dfrac{1}{4A^2B^2}\left[mA^2(\lambda^4\sin^4 r+ \lambda^4\sin^2 r+6\lambda^2+4)-8(n-1)B^2\lambda^4\sin^2 r \right]+1
\end{align*}
We show that the term in the above $[\cdot]$, denoted as I, is nonnegative for large $m$. We pick suitable positive terms to control the negative terms. When $r\in[\pi/4,\pi/2]$,
\begin{align*}
	\text{I}\ge &\ m\left( \lambda^8\sin^6r+6\lambda^6\sin^4r+4\lambda^4\sin^4r \right)\\
	&-8(n-1)\left( \lambda^8\sin^6r+2\lambda^6\sin^4r+\lambda^4\sin^2r \right)\ge 0
\end{align*}
provided $m\ge 8(n-1)$.
When $r\in [0,\pi/4]$, :
\begin{align*}
\text{I}\ge &\ m\left( \lambda^8\sin^6r+6\lambda^6\sin^4r+12\lambda^4\sin^2r\cos^2r \right)\\
&-8(n-1)\left( \lambda^8\sin^6r+2\lambda^6\sin^4r+\lambda^4\sin^2r \right)\ge 0
\end{align*}
provided $m\ge 8(n-1)$. This shows $Ric(H,H)\ge 1$.

Now we check $Ric(V,V)$. We have 
$$\dfrac{1-h'^2_\lambda}{h^2_\lambda}=\dfrac{A^3-\lambda^6\sin^2r}{\lambda^2A^2\cos^2r}\ge \dfrac{\lambda^6\sin^6r-\lambda^6\sin^2r}{\lambda^2A^2\cos^2r}=-\dfrac{\lambda^4\sin^2r(\sin^2r+1)}{A^2}.$$
\begin{align*}
	&Ric(V,V)= -\dfrac{h''_\lambda}{h_\lambda}+(n-2)\dfrac{1-h'^2_\lambda}{h^2_\lambda}-m\dfrac{f'_\lambda h'_\lambda}{f_\lambda h_\lambda}\\
	\ge&-\dfrac{2\lambda^4\sin^2r}{A^2}+1-(n-2)\dfrac{\lambda^4\sin^2r(\sin^2r+1)}{A^2}+m\dfrac{\lambda^2}{2A}\\
	\ge& \dfrac{\lambda^2}{A^2}\left[\dfrac{m}{2}\lambda^2\sin^2r-2n\lambda^2\sin^2r  \right]+1\ge1
\end{align*}
provided $m\ge 4n$. 

To summarize, we choose an integer $m\ge 8(n-1)$, then $Ric(g_\lambda)\ge 1$ holds for all $\lambda\ge 1$.

Let $\lambda\to\infty$. Then 
$$f_\lambda(r)\to 0,\quad h_\lambda(r)\to \tan^{-1} r.$$
Therefore, $(M,g_\lambda)$ Gromov-Hausdorff converges to the $n$-dimensional Grushin hemisphere as $\lambda\to\infty$. This proves Theorem A.

\begin{rem}
	The limit space has rectifiable dimension $n$ and Hausdorff dimension $2(n-1)$. Therefore, this also gives examples of Ricci limit spaces with Hausdorff dimension exceeding rectifiable dimension. When $n=2$, both Hausdorff and rectifiable dimension equal $2$, but the $2$-dimensional Hausdorff measure is not locally finite along the equator. 
\end{rem}

\begin{rem}
	Let $P$ be the north pole of the $n$-dimensional Grushin hemisphere. The closed metric ball $B_r(P)$ with radius $r$ is a convex subset, thus is RCD$(1,N)$ with a measure induced by the limit renormalized measure. As $r\to \pi/2$, $B_r(P)$ converges to the Grushin hemisphere. Note that $B_r(P)$ has Hausdorff dimension $n$ when $r<\pi/2$, and $2(n-1)$ when $r=\pi/2$. For $n\ge 3$, this shows that the Hausdorff dimension of RCD$(1,N)$ spaces does not have lower semi-continuity under the measured Gromov-Hausdorff convergence. This is different from the rectifiable dimension, which satisfies the lower semi-continuity \cite{Kit}.
\end{rem}

\begin{rem}
	Let $x$ be a point on the equator of the $n$-dimensional Grushin hemisphere. Then the tangent cone at $x$ is isometric to the $n$-dimensional Grushin halfspace with $\alpha=1$.
\end{rem}


\begin{thebibliography}{10}
	
\bibitem{Bel}	
A. Bella\"{\i}che.
\newblock The tangent space in sub-{R}iemannian geometry, In {\em Sub-Riemannian geometry}.
\newblock {\em Progr. Math.}, Vol. 144, Birkhäuser, Basel, 1-78, 1996.
	
\bibitem{BCGGJ}
U. Boscain, G. Charlot, J. P. Gauthier, S. Gu\'{e}rin, and H. R. Jayslin.
\newblock Optimal control in laser-induced population transfer for two-
and three-level quantum systems.
\newblock {\em J. Math. Phys.}, 43:2107--2132, 2002.
	
\bibitem{ChCo97}
J. Cheeger and T. H. Colding.
\newblock On the structure of spaces with {R}icci curvature bounded below. i.
\newblock {\em J. Differential Geom.}, 46(3):406--480, 1997.
	
\bibitem{DHPW}
X. Dai, S. Honda, J. Pan, and G. Wei.
\newblock Singular Weyl’s law with Ricci curvature bounded.
\newblock {\em arXiv}, 2208.13962, 2022.
	
\bibitem{Ju}
N. Juillet.
\newblock {Sub-Riemannian structures do not satisfy Riemannian Brunn–Minkowski inequalities}.
\newblock {\em Rev. Mat. Iberoam.}, 37(1):177-188, 2021.	
	
\bibitem{Kit} 
Y. Kitabeppu. 
\newblock A sufficient condition to a regular set being of positive measure on RCD spaces.
\newblock {\em Potential Anal.}, 51:179–196, 2019.	
	
\bibitem{MR}	
M. Magnabosco and T. Rossi.
\newblock Almost-Riemannian manifolds do not satisfy the CD condition.
\newblock {\em arXiv}, 2202.08775, 2022.
	
\bibitem{PanWei}
J. Pan and G. Wei.
\newblock Examples of Ricci limit spaces with non-integer Hausdorff dimension.
\newblock {\em Geom. Funct. Anal.}, 32, 676--685, 2022.	
		
		
\end{thebibliography}
\end{document}